\documentclass[11pt]{article} \title{\large\bf New Hermite-Hadamard Type Inequalities for Twice Differentiable Composite $(h-s)_2$-Convex  Functions}\author{Peter Olamide Olanipekun$^1$ and  Adesanmi Alao Mogbademu$^2$ }\date{}\usepackage{amssymb}\usepackage{amssymb}\usepackage[center]{titlesec}\usepackage[
  top=1.25in,
  bottom=1.25in,
  includeheadfoot,
  left=1.0in,
  right=1.0in,
  textwidth=6.6in,
  headsep=0.1in,
  headheight=0.2in]{geometry}\titlespacing*{\chapter}{0.00in}{0.01in}{2mm}\begin{document}\maketitle\noindent {\bf Abstract:}  In a recent paper [9], Ozdemir, Tunc and Akdemir defined two new classes of convex functions with which they proved some Hermite-Hadamard type inequalities. As an Open problem, they asked for conditions under which the composition of two functions belong to their newly defined class of convex functions and if Hermite-Hadamard type inequalities can be obtained. In this paper, we respond to the Open problems and  prove some new Hermite-Hadamard inequalities for twice differentiable composition whose second derivative is $((h-s)_{2}, I)$-convex. Our results are applied to some special means of real numbers.

\section{Introduction} The study of convex functions and some related inequalities has been an active area of research in mathematical ananlysis over the years. In an attempt to explain the concept of convex sets and convex functions, many researchers have discovered new classes of convex function (see for example [1], [4], [8], [6], [7], [11]) which in most cases, generalize the class of convex functions. Different inequalities which hold for the convex functions have been proven for other classes of convex functions. For several results in this direction, see [2],  [3], [5], [7] and the references therein.  We begin with the following definitions which are well known in literature.\\\\ {\bf Definition 1.1.} A function $f:I\subseteq \mathbb{R}\rightarrow\mathbb{R}$ is convex on an interval $I$ if for all $x, y\in I$ and $\alpha,\beta\in [0,1]$ such that $\alpha+\beta=1,$ the following inequality holds: $$f(\alpha x+\beta y)\leq \alpha f(x)+\beta f(y).$$ Let $f:I\subseteq \mathbb{R}\rightarrow\mathbb{R}$ be a convex function on the interval $I$ of real numbers and $a,b\in I$ with $a<b.$ The following double inequality is well known in literature as the Hermite-Hadamard inequality: $$f\left(\frac{a+b}{2}\right)\leq \frac{1}{b-a}\int_{a}^{b}f(x)\, dx\leq \frac{f(a)+f(b)}{2}.\eqno(1.1)$$ In 1978, Breckner introduced the class of $s-convex $ functions as follows: \\{\bf Definition 1.2.} A function $f:(0,\infty]\rightarrow [0,\infty]$ is $s$-convex in the second sense if for all $\alpha,\beta\in[0,1]$ such that $\alpha+\beta=1$, the following inequality holds $$f(\alpha x+\beta y)\leq \alpha^sf(x)+\beta^s f(y).$$ Another interesting class of convex functions (called the $h$-convex functions) was defined by Varosenic in [11] which generalizes many other classes, as follows: \\{\bf Definition 1.3.} Let $h:J\subseteq\mathbb{R}\rightarrow\mathbb{R}$ be a positive function, $h\not\equiv 0.$ A function $f:I\subseteq \mathbb{R}\rightarrow \mathbb{R}$ is $h$-convex if $f$ is nonnegative and for all $x,y\in I$, $\alpha,\beta\in [0,1]$ such that $\alpha+\beta=1$, we have $$f(\alpha x+\beta y)\leq h(\alpha)f(x)+h(\beta)f(y).$$ The following definitions are due to Ozdemir, Tunc and Akdemir [9] and are important for remaining part of this paper.\\\\{\bf Definition A.} Let $h:J\subset\mathbb{R}\rightarrow\mathbb{R}$ be a nonnegative function, $h\not\equiv 0.$ We say that $f:\mathbb{R}^+\cup\{0\}\rightarrow \mathbb{R}$ is an $(h-s)_1$-convex function in the first sense, or that $f$ belong to the class $SX((h-s)_1, I),$ if $f$ is non-negative and for all $x,y\in[0,\infty)=I$, $s\in(0, 1]$, $t\in[0,1]$, we have $$f(tx+(1-t)y)\leq h^s(t)f(x)+(1-h^s(t))f(y).\eqno(1.2)$$ If (1.2) is reversed, then $f$ is said to be $(h-s)_1$-concave function in the first sense, i.e $f\in SV((h-s)_1, I).$\\\\{\bf Definition B.} Let $h:J\subset \mathbb{R}\rightarrow\mathbb{R}$ be a nonnegative function, $h\not\equiv 0.$ We say that $f:\mathbb{R}^+\cup\{0\}\rightarrow\mathbb{R}$ is an $(h-s)_2$-convex function in the second sense, or that $f$ belong to the class $SX((h-s)_2, I),$ if $f$ is nonnegative and for all $u,v\in[0,\infty)=I$, $s\in(0,1]$, $t\in[0, 1]$ we have  $$f(tu+(1-t)v)\leq h^s(t)f(u)+h^s(1-t)f(v).\eqno(1.3)$$ If the inequality (1.3) is reversed, then $f$ is said to be $(h-s)_2$-concave function in the second sense, i.e., $f\in SV((h-s)_2, I).$   \section{Main Results}In [9], Ozdemir, Tunc and Akdemir defined two new classes of convex functions with which they proved some Hermite-Hadamard type inequalities. They also sought to know, as an Open problem, the conditions on the functions $f$ and $g$ for which the composition $f\circ g$ is $(h-s)_{1,2}$-convex on an interval $I$ and whether Hermite-Hadamard type Inequalities can be proven for the composition $f\circ g$?  In this section, we respond to the Open problems of Ozdemir, Tunc and Akdemir in [9].  We begin with the following observation.\\\\ (1). $SX(I)\subseteq SX((h-s)_{2}, I)$   whenever $h^s(t)\geq t$. \\(2). $K_{s}^{2}(I)\subseteq SX((h-s)_{2}, I)$  whenever $ h(t)\geq t$ \\(3). $SX(h,I)\subseteq SX((h-s)_{2}, I)$ whenever $h^{s-1}\geq 1$, $h>0$ \\where $SX(I), K_{s}^{2}(I)$ and $SX(h, I)$ denote the class of Convex, $s$-Convex, and $h$-Convex functions respectively. 
\\\\{\bf Theorem 2.1.} Let  $f\in SX((h-s)_{1}, I)$, $g:\mathbb{R}^+\cup\{0\}\rightarrow\mathbb{R}^+\cup\{0\}$ be a linear function, then the composition $f\circ g$ is $(h-s)_{1}$ convex on $I$.  
\\\\{\bf Proof.} \begin {eqnarray*}f\circ g(tx+(1-t)y)&=&f(tg(x)+(1-t)g(y))\\&\leq&h^s(t)f\circ g(x)+(1-h^s(t))f\circ g(y).\end{eqnarray*}  
\\\\{\bf Theorem 2.2.} Let $f$ be an increasing $(h-s)_{1}$ convex function on $I=[0,\infty)$. If $f$ is convex on $I$, then the composition $f\circ g$ is $(h-s)_{1}$ convex on the  interval $I$. 
\\\\{\bf Proof.} Since $g$ is convex implies \begin{eqnarray*}g(tx+(1-t)y)&\leq& tg(x)+(1-t)g(y) \\\mbox{Then}\hspace{1.4in} f\circ g(tx+(1-t)y)&\leq& f(tg(x)+(1-t)g(y))\\&\leq&h^s(t)f\circ g(x)+(1-h^s(t))f\circ g(y).\end{eqnarray*}
\\\\{\bf Theorem 2.3.} Let $f\in SX((h-s)_{2}, I)$, $g:\mathbb{R}^+\cup\{0\}\rightarrow\mathbb{R}^+\cup\{0\}$ be a linear function, then the composition $f\circ g$ is $(h-s)_{2}$ convex on $I$.
\\\\{\bf Proof.} \begin{eqnarray*}f\circ g(tu+(1-t)v)&=&f(tg(u)+(1-t)g(v))\\&\leq& h^s(t)f\circ g(u)+h^s(1-t)f\circ g(v). \end{eqnarray*}
\\\\{\bf Theorem 2.4.} Let $f$ be  an increasing $(h-s)_{2}$ convex function on $I=[0,\infty)$. If $g$ is convex on $I$, then the composition $f\circ g$ is $(h-s)_{2}$ convex on $I$. 
\\\\{\bf Proof.} \begin{eqnarray*}f\circ g(tu+(1-t)v)&\leq& f(tg(u)+(1-t)g(v))\\&\leq& h^s(t)f\circ g(u)+h^s(1-t)f\circ g(v)\end{eqnarray*}
\\\\{\bf Remark 2.5.} For the composition $f\circ g$ to belong to $SX((h-s)_{1,2}, I)$, at least $f$ must belong to $SX((h-s)_{1,2}, I)$. Suppose $f=g$ in Theorems 2.1- 2.4 above, then the self composition map $f^2$ is $(h-s)_{1,2}$ convex and by induction, $f^{k}, k=3, \cdots ,n$ is  $(h-s)_{1,2}$ convex. 
\\\\{\bf Theorem 2.6.} Let $f,g\in SX((h-s)_{2},I)$ such that $h^s(h^s(t))\leq h^s(t)$, then  the composition $f\circ g $ belongs to $SX((h-s)_{2}, I).$
\\\\{\bf Proof.} \begin{eqnarray*}f\circ g(tu+(1-t)v)&\leq& f\left(h^s(t)g(u)+h^s(1-t)g(v)\right)\\&\leq& h^s(h^s(t))f\circ g(u)+h^s(h^s(1-t))f\circ g(v) \\&\leq&h^s(t)f\circ g(u)+h^s(1-t)f\circ g(v).\end{eqnarray*} 
{\bf Theorem 2.7.} Let  $f,g\in SX((h-s)_{1},I)$ such that $h^s(h^s(t))= h^s(t)$, then  the composition $f\circ g $ belongs to $SX((h-s)_{1}, I).$
\\\\{\bf Proof.} \begin{eqnarray*}f\circ g(tx+(1-t)y)&\leq& f(h^s(t)g(x)+(1-h^s(t))g(y))\\&\leq&h^s(h^s(t))f\circ g(x)+(1-h^s(h^s(t)))f\circ g(y)\\&=& h^s(t)f\circ g(x)+(1-h^s(t))f\circ g(y).\end{eqnarray*}
{\bf Remark 2.8.}   Indeed,  $h(t)=\frac{1}{t}, s=\frac{1}{2}$ and $h(t)=t, s=1$ respectively satisfy the conditions on the function $h$  in Theorems 2.6 and 2.7. With the assumptions  given in Theorems 2.3, 2.4 and 2.6, we can prove Hermite-Hadamard inequality similar to those in literature for the composition $f\circ g$. In particular, the main theorems presented in [9] holds for the composition $f\circ g\in SX((h-s)_{2}, I)$. 
\\\\ In what follows, we denote the interior of $I$ by $I^\circ$,  the composition $f\circ g$ by $F$ and $\int_{0}^{1}h^s(t)$ by $K$.\\\\{\bf Theorem 2.9.}  Let $h\in L_{1}[0,1]$. If  $F:\mathbb{R}^+ \cup\{0\}\rightarrow \mathbb{R}^+\cup\{0\}$  is such that $F\in L_{1}[a, b],$ then under the assumptions of either Theorem 2.6 or 2.7, the following inequality holds:
$$\frac{1}{b-a}\int_{a}^{b} F(x)dx\leq K[F(a)+F(b)].\eqno(2.1)$$
\\\\{\bf Proof.} Clearly, from Theorems 2.6 and/or 2.7,  $F$ belongs to $SX((h-s)_{2}, I$. By setting $u=a$ and $v=b$ in (1.3), we obtain
$$F(ta+(1-t)b)\leq h^s(t)F(a)+h^s(1-t)F(b)\eqno(2.2)$$ Integrating both sides of (2.2) and setting $ta+(1-t)b=x$, we obtain 
$$\frac{1}{b-a}\int_{a}^{b} F(x)dx\leq F(a)K+ F(b)K.$$ This completes the proof.
\\\\{\bf Remark 2.10.} Observe that we have used the fact that $K=\int_{0}^{1}h^s(t)=\int_{0}^{1}h^s(1-t).$ 
\\\\{\bf Corollary 2.11.} Let $s=1$  in (2.1), then we obtain
$$\frac{1}{b-a}\int_{a}^{b}F(x)dx\leq [F(a)+F(b)]\int_{0}^{1}h(t)dt.$$
\\\\{\bf Corollary 2.12.} By choosing $h(t)=t$ in (2.1), we obtain the following \begin{eqnarray*}\frac{1}{b-a}\int_{a}^{b}F(x)dx&\leq&F(a)\int_{0}^{1}t^s dt+F(b)\int_{0}^{1}(1-t)^s dt\\&=&\frac{F(a)+F(b)}{s+1}\end{eqnarray*} which is the right hand side of inequality (11) in [9]. By choosing $s=1,$ we have the right hand side of  (1.1).
 \\\\Having established some conditions under which the composition of $f$ and $g$ can be $(h-s)_{1,2}$ convex on $I$, we now prove some Hermite-Hadamard type inequalities when the composition $f\circ g$ has an $(h-s)_{2}$ convex second derivative on $I$.
\\\\ We begin with the following lemmas.
\\\\{\bf Lemma 2.13.} Let $F:I^\circ\subseteq\mathbb{R}\rightarrow \mathbb{R}$ be a twice differentiable mapping on $I^\circ$, where $a, b\in I^\circ$ with $a<b$. If $F''\in L[a, b],$ then the following inequality holds $$\frac{F'(a)+F'(b)}{2}-\frac{1}{b-a}\int_{a}^{b}F'(x)dx=\frac{b-a}{2}\int_{0}^{1}(1-2t)F''(ta+(1-t)b)dt.\eqno(2.3)$$ 
\\\\{\bf Proof.} Using integration by part and the substitution $x=ta+(1-t)b$, where $t\in [0,1], $ we have  \begin{eqnarray*}\int_{0}^{1} (1-2t)F''(ta+(1-t)b) dt &=& (1-2t)\left.\frac{F'(ta+(1-t)b)}{a-b}\right|_{0}^{1}+ 2\int_{0}^{1}\frac{F'(ta+(1-t)b)}{a-b} dt\\&=& \frac{F'(a)+F'(b)}{b-a}-\frac{2}{(b-a)^2}\int_{a}^{b} F'(x) dx. \end{eqnarray*} By re-arranging we obtain (2.3).
\\\\{\bf Lemma 2.14.} Let $F:I^\circ\subseteq\mathbb{R}\rightarrow \mathbb{R}$  be a twice differentiable mapping on $I^\circ$ where $a, b\in I^\circ$ with $a<b$. If $F''\in L[a, b]$, then $$\frac{F(a)+F(b)}{2}-\frac{1}{b-a}\int_{a}^{b} F(x) dx= \frac{(b-a)^2}{2}\int_{0}^{1}(t-t^2)F''(ta+(1-t)b)dt.\eqno(2.4)$$
\\\\{\bf Proof.} By applying integration by parts twice and the substitution $x=ta+(1-t)b$, we have \begin{eqnarray*}\int_{0}^{1} (t-t^2)F''(ta+(1-t)b) dt&=& (t-t^2)\left.\frac{F'(ta+(1-t)b)}{a-b}\right|_{0}^{1}-\int_{0}^{1}\frac{F'(ta+(1-t)b)}{a-b} (1-2t) dt \\&=&\frac{F(a)+F(b)}{(b-a)^2}-\frac{2}{(b-a)^3}\int_{a}^{b}F(x)dx.\end{eqnarray*} By re-arranging, we obtain (2.4).
\\\\{\bf Theorem 2.15.} Let $F:I^\circ\subseteq\mathbb{R}\rightarrow \mathbb{R}$ be a twice differentiable mapping on $I^\circ$, where $a,b\in I^\circ$ with $a<b$. If $|F''|$ is $(h-s)_{2}$ convex on $[a, b]$, then \begin{eqnarray*}\left|\frac{F'(a)+F'(b)}{2}-\frac{1}{b-a}\int_{a}^{b}F'(x)dx\right|&\leq& \frac{b-a}{2}\left(F''(a)+F''(b)\right)\int_{0}^{1}\left|1-2t\right|h^s(t)dt\\&\leq&\frac{b-a}{2}\left(F''(a)+F''(b)\right) \left(K+2\int_{0}^{1}t\, h^s(t)dt\right).\end{eqnarray*} 
{\bf Proof.} Using Lemma 2.13, \begin{eqnarray*}\left|\frac{F'(a)+F'(b)}{2}-\frac{1}{b-a}\int_{a}^{b}F'(x)dx\right|&=&\left|\frac{b-a}{2}\int_{0}^{1}(1-2t)F''(ta+(1-t)b)dt\right|\\&\leq&\frac{b-a}{2}\int_{0}^{1}|1-2t|\left|F''(ta+(1-t)b)\right|dt\\&\leq&\frac{b-a}{2}\int_{0}^{1}|1-2t|\left(h^s(t)\left|F''(a)\right|+h^s(1-t)\left|F''(b)\right|\right)dt\\&\leq&\frac{b-a}{2}\int_{0}^{1}(1+2t)\left(h^s(t)\left|F''(a)\right|+h^s(1-t)\left|F''(b)\right|\right)dt\\&=&\frac{b-a}{2}\left(F''(a)+F''(b)\right) \left(K+2\int_{0}^{1}t\, h^s(t)dt\right).\end{eqnarray*} This completes the proof.
\\\\{\bf Corollary 2.16.} Let $F:I^\circ\subseteq\mathbb{R}\rightarrow \mathbb{R}$ be a twice differentiable mapping on $I^\circ$, where $a,b\in I^\circ$ with $a<b$. If $|F''|$ is $(h-s)_{2}$ convex on $[a, b]$, such that $s=1, h(t)=t$, then $$\left|\frac{F'(a)+F'(b)}{2}-\frac{1}{b-a}\int_{a}^{b}F'(x)dx\right|\leq\frac{(b-a)(|F'(a)|+|F'(b)|)}{8}.$$
{\bf Theorem 2.17.} Let    $F:I^\circ\subseteq\mathbb{R}\rightarrow \mathbb{R}$ be a twice differentiable mapping on $I^\circ$, where $a,b\in I^\circ$ with $a<b$ and let $p>1$. If the  mapping $|F''|^\frac{p}{p-1}$ is $(h-s)_{2}$ convex on $[a, b]$ then $$\left|\frac{F'(a)+F'(b)}{2}-\frac{1}{b-a}\int_{a}^{b}F'(x)dx\right|\leq \frac{(b-a)}{2(p+1)^{\frac{1}{p}}}\left[K\left(\left|F''(a)\right|^\frac{p}{p-1}+\left|F''(b)\right|^\frac{p}{p-1}\right)\right]^\frac{p-1}{p}.\eqno(2.5)$$ 
{\bf Proof.} Using Lemma 2.13 and the Holder's inequality, we have 
\begin{eqnarray*}
\left|
\frac{F'(a)+F'(b)}{2}-\frac{1}{b-a}\int_{a}^{b}F'(x)dx\right|&\leq& \frac{b-a}{2}\int_{0}^{1}|1-2t||F''(ta+(1-t)b)| \,dt\\&\leq&\left(\frac{b-a}{2}\right)\left(\int_{0}^{1}\left|1-2t\right|^p\, dt\right)^\frac{1}{p}\left(\int_{0}^{1}\left|F''(ta+(1-t)b)\right|^\frac{p}{p-1}\, dt\right)^\frac{p-1}{p}\hspace{0.1in}(2.6)
\end{eqnarray*} By the   $(h-s)_{2}$ convexity of $|F''|^\frac{p}{p-1}$ we have,
\begin{eqnarray*}
\int_{0}^{1}\left|F''(ta+(1-t)b)\right|^\frac{p}{p-1}\, dt&\leq& \int_{0}^{1}\left(h^s(t)\left|F''(a)\right|^\frac{p}{p-1}+h^s(1-t)\left|F''(b)\right|^\frac{p}{p-1}\right)\, dt \\&=& \left[\left|F''(a)\right|^\frac{p}{p-1}+\left|F''(b)\right|^\frac{p}{p-1}\right]\int_{0}^{1}h^s(t)\, dt \hspace{1.36in} (2.7)
\end{eqnarray*} Also, since $$\int_{0}^{1}|1-2t|^p \, dt=\int_{0}^{\frac{1}{2}}(1-2t)^p+\int_{\frac{1}{2}}^{1}(2t-1)^p\,dt=2\int_{0}^{\frac{1}{2}}(1-2t)p\, dt=\frac{1}{p+1}\eqno(2.8)$$ By using (2.7) and (2.8) in (2.6) we obtain (2.5).
\\\\{\bf Corollary 2.18.} Let $p=2,$ and $K=K(s,t)=K(1,t)=t^2$, then 
\begin{eqnarray*}
\left|\frac{F'(a)+F'(b)}{2}-\frac{1}{b-a}\int_{a}^{b}F'(x)\, dx\right|&\leq& \frac{(b-a)\sqrt{3}}{6}\left[t^2\left(F''(a))^2+(F'(b))^2\right)\right]^{\frac{1}{2}}\\&=&\frac{(b-a)\sqrt{3}t}{6}\left(\left(F''(a)+F''(b)\right)^2-2F''(a)''F''(b)\right)^{\frac{1}{2}}.
\end{eqnarray*}
\section{Applications to Some Special Means} The following special means are well known in literature (See for example [5], [10]). \\  Arithmetic Mean $$A(a, b)=\frac{a+b}{2}.$$Harmonic Mean $$H(a, b)=\frac{2ab}{a+b}.$$ Logarithmic Mean $$L(a, b)=\frac{b-a}{ \ln b- \ln a}$$ Power Mean of order $p$ $$A_{p}(a, b)=\left(\frac{a^p+b^p}{2}\right)^{\frac{1}{p}}$$

As an application of our results, we present the following Propositions.\\\\{\bf Proposition 3.1.} Let $a, b\in \mathbb{R}$, $a<b$ with $p>1$. Then the following inequality holds: $$\left|\frac{1}{L(a, b)}-\frac{1}{H(a, b)}\right|\leq \frac{b-a}{2^{\frac{3p-2}{p}}(p+1)^{\frac{1}{p}}}A_{\frac{2p}{1-p}}^{2}(a, b).$$ {\bf Proof.} It follows immediately from Theorem 2.17 when $F(x)=-\ln x$, $x>0$ and $h(t)=t.$\\\\{\bf Proposition 3.2.} Let $a, b\in \mathbb{R}$, $a< b$. Then  $$\left|\frac{1}{L(a, b)}-\frac{1}{H(a, b)}\right|\leq\frac{b^2-a^2}{8ab}.$${\bf Proof.} It follows immediately from Corollary 2.16 by using $F(x)=-\ln x$, $x>0$.\\\\{\bf Proposition 2.3.} Let $a,b\in\mathbb{R}$, $a<b,$ and $n\in\mathbb{N}$, $n\geq 2.$ Then for any $p>1$, the following holds: $$\left|nA\left(a^{n-1},b^{n-1}\right)-L_{n}^{n}(a,b)\right|\leq \frac{b-a}{2(p+1)^{\frac{1}{p}}}(n^2-n)^{\frac{p-1}{p}}\left[A\left(\left|a\right|^{(n-2)p/p-1},\left|b\right|^{(n-2)p/p-1}\right)\right]^{\frac{p-1}{p}}.$$ {\bf Proof.} It follows from Theorem 2.17 by using $F(x)=x^n$, $x\in \mathbb{R}$, $n\geq2.$ 
 \bibliographystyle{amsplain}{\footnotesize
Research Group in Mathematics and Applications\\Department of Mathematics\\
University of Lagos, Lagos, Nigeria\\[-1mm]
e-mail: {\tt $^{1}$polanipekun@yahoo.com, $^{2}$amogbademu@unilag.edu.ng } } \\ [2mm]

\end{document}